\magnification=1100
\input amstex
\documentstyle{amsppt}
\def\leftitem#1{\item{\hbox to\parindent{\enspace#1\hfill}}}

\def\dim{\operatorname{dim}}

\def\proclaim#1{\par\vskip.25cm\noindent{\bf#1. }\begingroup\it}
\def\endproclaim{\endgroup\par\vskip.25cm}
\leftheadtext{Shen, Wang} \rightheadtext{almost inner derivations of Lie algebras}

\topmatter
\title Almost Inner Derivations of Affinizations of Minimal Q-graded Subalgebras
\endtitle
\author  Yaxin Shen$^{\dagger}$, Xiandong Wang$^{\ddag}$
\endauthor
\affil   College of Mathematics and Statistics\\
       Qingdao University, Qingdao 266071\\
       P. R. China\\
       $^{\dag}$Email:  shenyaxin\@qdu.edu.cn \,\,\,\,\,\,\,\,\,\,\\
       $^{\ddag}$Email: wangxiandong\@qdu.edu.cn\\
\endaffil
\thanks Research supported by NSF of China (No.11472144).
\endthanks
\keywords derivations, almost inner derivations, minimal Q-graded subalgebras, loop algebras, affinizations
 \endkeywords
\subjclass\nofrills 2020 {\it Mathematics Subject Classification.}
 17B40, 17B65\endsubjclass
 \abstract Minimal Q-graded subalgebras of semisimple Lie algebras are introduced,
 and it is proved that their derived algebras are abelian.
 Almost inner derivations of minimal Q-graded subalgebras are investigated, they are all inner derivations.
 Based on these Lie algebras, a decomposition formula is obtained for derivations of loop algebras,
 and almost inner derivations of affinizations are determined.
 \endabstract
\endtopmatter
\document

\smallskip\bigskip
\subhead 1. \ \ Introduction
\endsubhead
\medskip

Almost inner derivations of Lie algebras originate from almost inner automorphisms of the corresponding Lie groups,
the latter is used to resolve isospectral problems of compact Riemannian manifolds.
The early invesgators for these related topics include J.Milnor, C.Gordan and E.Wilson, see [9,14] for their discussions
on basic results and some concrete examples.

Recently, D.Burde, K.Dekimpe and B.Verbeke began a systematic study of almost inner derivations of Lie algebras,
they formulated some fundamental features and gave many interesting examples[6-8]. Their methods are adequate to determine
whether or not the space of almost inner derivations equal to the space of inner derivations for many (mainly nilpotent) Lie algebras.
From these results one gets some impressions of the picture of this anticipating subject.

In the present paper, we provide new classes of Lie algebras in which almost inner derivations are all determined.
Compared with that in the above literature, the minimal Q-graded subalgebras are solvable but not nilpotent,
and they are some kinds of modified Q-graded Lie algebras as introduced by I.Kenji and K.Yoshiyuki in [12].
In section 2, we give the definition of a minimal Q-graded subalgebra of a semisimple Lie algebra,
and prove that its derived algebra is commutative, this is essential for the main results of the paper.
In section 3, we first recall some basic definitions on almost inner derivations, then prove two key lemmas,
and finally determine all the almost inner derivations for minimal Q-graded subalgebras.
In section 4, we consider the loop algebra $L\otimes S$ \,($S=\Bbb{F}[t,t^{-1}]$) of a minimal Q-graded subalgebra $L$,
and describe derivations and almost inner derivations for
this infinite dimensional Lie algebra using information on tensor product of \,(non-associative) algebras,
on the centroid of minimal Q-graded subalgebras, etc. We thus obtain a familiar decomposition
$$ Der(L\otimes S)=Der(L)\otimes S\oplus Cent(L)\otimes Der(S) $$
that usually holds under the assumption that the algebra $L$ is perfect.
Finally, in section 5, all almost inner derivations of the
affinization $\tilde{L}=L\otimes S\oplus \Bbb{F}K$ \,(usually denoted by $AID(\tilde L)$)
of a minimal Q-graded subalgebra $L$ with $\dim L=2\dim[L,L]$ are determined:

The quotient vector space $AID(\tilde{L})/Inn(\tilde{L})$ consists of elements of the form
$$ \sum_{1\leq i\leq l, j\in \Bbb{Z}}a_{ij}D_{ij}+Inn(\tilde{L}),\,\, a_{ij}\in \Bbb{F}, $$
where the derivation $D_{ij}$ is almost inner and determined by the following rules:
$$ D_{ij}([L,L]\otimes S)=0, \,\, D_{ij}(K)=0, $$
$$ D_{ij}(h_m\otimes t^n)=\delta_{im}\delta_{jn}K,\,\, 1\leq m\leq l, \,\, n\in \Bbb{Z},$$
here $h_1,\cdots ,h_l$ is a basis of the maximal torus $H$ which is contained in the subalgebra $L$.
Moreover, the infinite subset $\{D_{ij}+Inn(\tilde{L}); 1\leq i\leq l, j\in \Bbb{Z}\}$ of the quotient space
$ AID(\tilde{L})/Inn(\tilde{L}) $ is linearly independent.

 \vskip 5pt

\smallskip\bigskip
\subhead 2. \ \ Minimal Q-graded subalgebras of semisimple Lie algebras
\endsubhead
\medskip

First we fix some notations: Let $\Bbb{F}$ be an algebraically closed field with characteristic $0$,
$\Bbb{L}$ a finite dimensional semisimple Lie algebra over $\Bbb{F}$.
Also denote by $H$ a maximal torus, $\Phi$ the corresponding root system, and $\triangle=\{\alpha_{1},\alpha_{2},\cdots,\alpha_{l}\}$ a basis of $\Phi$,
hence we have the root space decomposition of $\Bbb{L}$:
$$ \Bbb{L}=H\oplus \bigoplus_{\alpha\in \Phi}L_{\alpha},$$
where $L_{\alpha}$ is the root space of $\Bbb{L}$ with respect to the root $\alpha$. For more terminologies and elementary results on Lie algebras,
take [10-11] as basic references.

\proclaim{Definition 2.1}
Given a semisimple Lie algebra $\Bbb{L}$, a subalgebra $L$ of $\Bbb{L}$ including the maximal torus $H$ has the following form:
$$ L=H\oplus\bigoplus_{\alpha\in{\Psi}}L_{\alpha}, \,\,\, \Psi{\subseteq}{\Phi}.$$
If the set $\Psi$ is a spanning subset of $Q=\Bbb{Z}{\Phi}$ as a free abelian group in the sense that $Q=span_{\Bbb{Z}}{\Psi}$,
then $L$ is called a $Q$-graded subalgebra of $\Bbb{L}$.

Moreover, if there is no $Q$-graded subalgebra $L_{1}$ of $\Bbb{L}$ such that ${L_{1}}\subsetneq{L}$,
then $L$ is called a minimal $Q$-graded subalgebra of $\Bbb{L}$.
\endproclaim

Remark: The above definition for a Q-graded subalgebra is motivated from the concept of Q-graded Lie algebras in [12].
Notice also that these Lie algebras are not commonly referred in literature as root graded Lie algebras,
which are discussed intensively by many authors, see [3,5,13].
The main difference is that root graded Lie algebras are perfect, but minimal Q-graded ones are not.
Perfect algebras have automatically commutative centroid, but for minimal Q-graded subalgebras it needs some effort to
calculate its centroid, and this will be done in the next section.

The following two examples exhibit some concrete minimal Q-graded subalgebras in which the equality $\dim H=\dim [L,L]$ holds,
and the inequality $\dim H\leq \dim [L,L]$ holds in general by the very definition.

\proclaim{Example 2.2}
Let $\Bbb{L}$ be the classical $B_{2}$-type simple Lie algebra with root space decomposition: $\Bbb{L}=H\oplus\oplus_{\alpha\in{\Phi}}L_{\alpha}$,
where the root system $\Phi$ is as follows
$$\Phi={\{\alpha,\beta,\alpha+\beta,2\alpha+\beta,-\alpha,-\beta,-\alpha-\beta,-2\alpha-\beta\}},$$
for which we choose a basis $\Delta=\{\alpha,\beta\}\subseteq\Phi$.
It can be checked easily that $\Bbb{L}$ has and only has the following minimal $Q$-graded subalgebras:

$H\oplus{L_{\alpha}}\oplus{L_{2\alpha+\beta}}$; \,\,\,\,\,\,\,\,\,\,\,\, $H\oplus{L_{\alpha}}\oplus{L_{-\beta}}$;
\,\,\,\,\,\,\,\,\,\,\,\,\,\,\,\,\,\,\, $H\oplus{L_{\beta}}\oplus{L_{\alpha+\beta}}$;

$H\oplus{L_{\beta}}\oplus{L_{-\alpha}}$;\,\,\,\,\,\,\,\,\,\,\,\,\,\,\,\,\,\,\, $H\oplus{L_{\alpha+\beta}}\oplus{L_{2\alpha+\beta}}$;
\,\,\,\,\,\,\,\, $H\oplus{L_{-\alpha}}\oplus{L_{-2\alpha-\beta}}$;

$H\oplus{L_{-\beta}}\oplus{L_{-\alpha-\beta}}$;\,\,\,\,\,\,\,\,$H\oplus{L_{-\alpha-\beta}}\oplus{L_{-2\alpha-\beta}}$.
\endproclaim

\proclaim{Example 2.3}
Let $\Bbb{L}$ be the classical $A_{l}$-type simple Lie algebra, it has root space decomposition: $\Bbb{L}=H\oplus\oplus_{\alpha\in{\Phi}}L_{\alpha}$.
By the well-known functions $\varepsilon_i: H\rightarrow \Bbb{F}$,
a basis $\triangle=\{\alpha_1,\alpha_2,\cdots,\alpha_l\}$ of $\Phi$ is defined naturally as follows:
$$\Phi=\{\varepsilon_{i}-\varepsilon_{j}; 1\leq{i}\neq{j}\leq{l+1}\},$$
$$\alpha_{i}=\varepsilon_{i}-\varepsilon_{i+1}, 1\leq{i}\leq{l},\,\,\, \,\,\,\,$$
$$\varepsilon_{i}-\varepsilon_{j}=\varepsilon_{i}-\varepsilon_{i+1}+\cdots+\varepsilon_{j-1}-\varepsilon_{j},$$
$$\,\,\, =\alpha_{i}+\cdots+\alpha_{j},\,\,\,i<j.$$
From this we see that the positive roots are: $\varepsilon_{i}-\varepsilon_{j}, i<j$,
and we have some minimal $Q$-graded subalgebras of $\Bbb{L}$ below
$$ H\oplus L_{\alpha_1}\oplus L_{-\alpha_2}\oplus \cdots \oplus L_{(-1)^{l-1}\alpha_l},$$
$$H\oplus{L_{\alpha_{1}}}\oplus{L_{\alpha_{1}+\alpha_{2}}}\oplus\cdots\oplus{L_{\alpha_{1}+\cdots+\alpha_{l}}},$$
$$H\oplus{L_{\alpha_{2}}}\oplus{L_{\alpha_{1}+\alpha_{2}}}\oplus\cdots\oplus{L_{\alpha_{1}+\cdots+\alpha_{l}}},$$
$$H\oplus{L_{\alpha_{1}}}\oplus{L_{\alpha_{3}}}\oplus{L_{\alpha_{1}+\alpha_{2}+\alpha_{3}}}\oplus\cdots\oplus{L_{\alpha_{1}+\cdots+\alpha_{l}}},$$
$$H\oplus{L_{\alpha_{2}}}\oplus{L_{\alpha_2+\alpha_{3}}}\oplus{L_{\alpha_{1}+\alpha_{2}+\alpha_{3}}}\oplus\cdots\oplus{L_{\alpha_{1}+\cdots+\alpha_{l}}}.$$
\endproclaim

Now we can state and prove the main result of this section, which says roughly that a minimal Q-graded subalgebra is ``almost abelian".
\proclaim{Theorem 2.4}
Let $\Bbb{L}$ be a semisimple Lie algebra, and $L$ any minimal $Q$-graded subalgebra: $L=H\oplus\oplus_{\alpha\in{\Psi}}L_{\alpha}$, $\Psi{\subseteq}{\Phi}$.
Set $I=[L,L]=\oplus_{\alpha\in{\Psi}}L_{\alpha}$, then $I$ is abelian.
In this case, the Lie algebra $L$ is also called metabelian.
\endproclaim

Before proving this theorem, let's first recall the height of a root:
For any given root $\beta\in\Phi$, its height with respect to the basis $\Delta=\{\alpha_{1},\alpha_{2},\cdots,\alpha_{l}\}$
is defined separately below

(1) If $\beta\in\Phi^{+}$, the height of $\beta=\Sigma_{i=1}^{l}k_{i}\alpha_i$ is $ht(\beta)=\Sigma_{i=1}^lk_{i}$,

(2) If $\beta\in\Phi^{-}$, just define $ht(\beta)=-ht(\beta)$.

\demo{Proof} We divide the whole proof into three parts, and reorder the roots in $\Psi$ according to their heights if necessary.

1) $\Psi=\{\beta_1,\beta_2,\cdots,\beta_m\}\subseteq\Phi^{+}$, $ht(\beta_{1})\leq ht(\beta_{2})\leq\cdots\leq ht(\beta_{m})$.

a) If $[L_{\beta_1},L_{\beta_2}]\neq0$, then there exists some $k\geq3$ satisfying $[L_{\beta_1},L_{\beta_2}]=L_{\beta_{k}}$, and so $\beta_{1}+\beta_{2}=\beta_{k}$. According to the assumption on heights, $\beta_{i}+\beta_{j}\neq\beta_{1}$, $\forall{i,j\geq1}$.
Set $\Psi_{1}=\{\beta_2,\cdots,\beta_m\}$, then $H\oplus\oplus_{\alpha\in{\Psi_{1}}}L_{\alpha}$ is a subalgebra and obviously it is also Q-graded.
Similarly, if $[L_{\beta_{1}},L_{\beta_{j}}]\neq0$ for some $j\geq2$, then $H\oplus\oplus_{\alpha\in{\Psi_1}}L_{\alpha}$ is a $Q$-graded subalgebra.

b) If $[L_{\beta_{1}},L_{\beta_{j}}]=0$ for any $j\geq2$ and $[L_{\beta_{2}},L_{\beta_{j}}]\neq 0$ for some $j\geq3$,
then we must have $\beta_{2}+\beta_{j}\in{\Psi}$.
Set $\Psi_{1}=\Psi\backslash\{\beta_{2}\}$, and notice that in this case, $\beta_{1}+\beta_{j}$ is not a root for any $j\geq2$
and also $\beta_{i}+\beta_{j}\neq\beta_{2}$, $\forall{i,j\geq2}$, hence $H\oplus\oplus_{\alpha\in{\Psi_{1}}}L_{\alpha}$ is a $Q$-graded subalgebra.

c) If $[L_{\beta_{i}},L_{\beta_{j}}]=0$, $1\leq{i}\leq{r-1}$, $\forall{j}\geq1$ and $[L_{\beta_{r}},L_{\beta_{j}}]\neq0$, $\exists{j}\geq{r+1}$,
we similarly define
 $\Psi_{1}=\Psi\backslash\{\beta_{r}\}$. In this case, $\beta_{i}+\beta_{j}$  is not a root, $1\leq{i}\leq{r-1}$, $\forall{j}\geq1$, and  $H\oplus\oplus_{\alpha\in{\Psi_{1}}}L_{\alpha}$ is a $Q$-graded subalgebra.

From a),b),c) above, we get a contradiction to the minimality of the Q-graded subalgebra $L$,
therefore, $[L_{\beta_{i}},L_{\beta_{j}}]=0$, $\forall{i,j\geq1}$, and $\oplus_{\alpha\in{\Psi}}L_{\alpha}$ is abelian.

2) $\Psi=\{\beta_{1},\beta_{2},\cdots,\beta_{m}\}\subseteq\Phi^{-}$, this is the case analogous to that in 1),
and we can follow the same discussions to obtain that $\oplus_{\alpha\in{\Psi}}L_{\alpha}$ is abelian.

3) Finally, suppose the subset $\Psi$ contains both positive and negative roots.
In this case, we first give a specific decomposition of $L$ as follows:
$$L=H\oplus{L_{{\beta'}_{s'}}}\oplus\cdots\oplus{L_{{\beta'}_1}}\oplus{L_{{-\gamma'}_1}}\oplus\cdots\oplus{L_{{-\gamma'}_{r'}}},$$
here $\{{\beta'}_{s'}, \cdots ,{\beta'}_1, {\gamma'}_1,\cdots, {\gamma'}_{r'}\}\subseteq \Phi^+$, and $s'\geq 1, r'\geq 1$.

We will normalize the above decomposition by deleting some root spaces of $L$ according to the principals:
If some root ${\gamma'}_{i'}$ can be expressed as
an integral linear combination of ${{\beta'}_{1}},\cdots,{{\beta'}_{s'}}$, then delete the corresponding root space in $L$,
and denote by ${\gamma}_{1}$ the root in $\{{\gamma'}_1,\cdots, {\gamma'}_{r'}\}$ which cannot be expressed as an
integral linear combination of ${{\beta'}_{1}},\cdots,{{\beta'}_{s'}}$ with maximal height.

Similarly, if some root ${\gamma'}_{j'}$ can be expressed as an integral linear combination of
${{\gamma}_{1}},{{\beta'}_{1}},\cdots,{{\beta'}_{s'}}$, then delete the corresponding root space in $L$,
and denote by ${\gamma}_{2}$ the root in  $\{{\gamma'}_1,\cdots, {\gamma'}_{r'}\}$ which cannot be expressed as an
integral linear combination of $\gamma_1,{{\beta'}_{1}},\cdots,{{\beta'}_{s'}}$ with maximal height.

Proceeding in this way, the Q-graded subalgebra $L$ becomes to the subspace
$$ H\oplus{L_{{\beta'}_{s'}}}\oplus\cdots\oplus{L_{{\beta'}_{1}}}\oplus{L_{{-\gamma}_{1}}}\oplus\cdots\oplus{L_{{-\gamma}_{r}}}.$$
This subspace will be normalized further by a similar consideration as described previously,
and we get a subspace like this:
$$ H\oplus{L_{{\beta}_{s}}}\oplus\cdots\oplus{L_{{\beta}_{1}}}\oplus{L_{{-\gamma}_{1}}}\oplus\cdots\oplus{L_{{-\gamma}_{r}}},$$
$$ ht(\beta_{1})\geq ht(\beta_{2})\geq\cdots\geq ht(\beta_{s}),\,\, ht(\gamma_{1})\geq ht(\gamma_{2})\geq\cdots\geq ht(\gamma_{r}).$$

Now it suffices to prove that the subspace ${L_{{\beta}_{s}}}\oplus\cdots\oplus{L_{{\beta}_{1}}}\oplus{L_{{-\gamma}_{1}}}\oplus\cdots\oplus{L_{{-\gamma}_{r}}}$
is actually an abelian subalgebra. If so, it must be the original $[L,L]$ by minimality, as the Q-gradation requirement is obviously valid.

a) Claim: The subspace ${L_{{\beta}_{s}}}\oplus\cdots\oplus{L_{{\beta}_{1}}}$ is an abelian subalgebra.

i) If for some $j\geq 2$ and $\beta'\in \Psi$, we have $[L_{\beta_{1}},L_{\beta_{j}}]=L_{\beta'}$, then $\beta_{1}+\beta_{j}={\beta'}$.
According to the assumption about the height,
${\beta'}$ can be expressed as an integral linear combination of ${{\gamma}_{1}},\cdots,{{\gamma}_{r}}$,
therefore ${\beta_{j}}$ can be expressed as an integral linear combination of ${\beta_{1}},{{\gamma}_{1}},\cdots,{{\gamma}_{r}}$,
and this is a contradiction. Hence $[L_{\beta_{1}},L_{\beta_{j}}]=0$ for any $j\geq2$.

ii) If for some $j\geq 3$ and $\beta'\in \Psi$, we have $[L_{\beta_{2}},L_{\beta_{j}}]=L_{\beta'}, \beta_{2}+\beta_{j}={\beta'}$,
then by similar discussions as in i), we will derive a contradiction, hence $[L_{\beta_{2}},L_{\beta_{j}}]=0$ for all $j\geq3$.

Repeating the above processes, it can be deduced that the subspace ${L_{{\beta}_{s}}}\oplus\cdots\oplus{L_{{\beta}_{1}}}$
is indeed an abelian subalgebra as desired.

b) Claim: The subspace ${L_{-{\gamma}_{1}}}\oplus\cdots\oplus{L_{-{\gamma}_{r}}}$ is an abelian subalgebra.

The proof of this claim is omitted, since it can be checked by a similar discussion as that used in a).

c) Claim: $[L_{\beta_{i}},L_{-{\gamma}_{j}}]=0$, $1\leq{i}\leq{s}$, $1\leq{j}\leq{r}$. Contrary to this claim,
suppose $[L_{\beta_{i_0}},L_{-{\gamma}_{j_0}}]\neq0$ for some $i_0,j_0$, we will get a contradiction
from the following calculations in both cases:

i) If $[L_{\beta_{i_0}},L_{-{\gamma}_{j_0}}]=L_{\beta'}$ for some $\beta'\in \Psi$, and so ${\beta_{i_0}}-{\gamma_{j_0}}={\beta'}$,
then the root ${\gamma}_{j_0}$ can be expressed as an integer linear combination of ${{\beta'}_{1}},\cdots,{{\beta'}_{s'}}$.
This contradicts the hypothesis that ${\gamma}_{j_0}$ cannot be expressed as an integer linear combination of roots ${{\gamma}_{1}},\cdots,{{\gamma}_{j_0-1}},{{\beta'}_{1}},\cdots,{{\beta'}_{s'}}$.

ii) If $[L_{\beta_{i_0}},L_{-{\gamma}_{j_0}}]=L_{-{\gamma'}}$ for some $\gamma'\in \Psi$, and so ${\beta_{i_0}}-{\gamma_{j_0}}={-{\gamma'}}$,
then the root ${\beta}_{i_0}$ can be expressed as an integer linear combination of ${{\gamma'}_{1}},\cdots,{{\gamma'}_{r'}}$.
This contradicts the hypothesis that ${\beta}_{i_0}$ cannot be expressed as an integer linear combination of roots ${\beta_{1}},\cdots,{\beta}_{i_0-1},{\gamma'}_{1},\cdots,{\gamma'}_{r'}$.

From the three claims above, we deduce that the subspace
${L_{{\beta}_{s}}}\oplus\cdots\oplus{L_{{\beta}_{1}}}\oplus{L_{{-\gamma}_{1}}}\oplus\cdots\oplus{L_{{-\gamma}_{r}}}$ is an abelian subalgebra,
it is Q-graded, and equal to $[L,L]$ by minimality condition.  The proof of the theorem is finished.
\hfill $ $ \qed
\enddemo

 \vskip 5pt

\smallskip\bigskip
\subhead 3. \ \  Almost inner derivations of Lie algebras
\endsubhead
\medskip

A systematic discussion of almost inner derivations of Lie algebras over arbitrary fields is carried out in [6-8],
and we first recall the basic definitions and properties which will be used in this paper.

\proclaim{Definition 3.1}
A derivation $D$ of a Lie algebra $L$ is called an almost inner derivation if $D(x)\in[x,L], \forall{x}\in{L}$.
The vector space consisting of all almost inner derivations of $L$ is denoted by $AID(L)$.

A derivation is almost inner if and only if it coincides on each one-dimensional subspace with an inner derivation.
In particular, the set $Inn(L)$ of all inner derivations of $L$ is a subspace of $AID(L)$.
\endproclaim

\proclaim{Proposition 3.2}
We have the following inclusions between Lie subalgebras of the general linear Lie algebra $gl(L)=End(L)$:
$${Inn(L)}\subseteq{AID(L)}\subseteq{Der(L)},$$
here $Der(L)$ is the Lie algebra consisting of all derivations of the Lie algebra $L$.
\endproclaim
\demo{Proof} See [6] for a direct proof of this basic property.
\hfill $ $ \qed
\enddemo

Our goal in this section is to determine all almost inner derivations of minimal Q-graded subalgebras of semisimple Lie algebras,
we first prove the following two key lemmas for this purpose.

\proclaim{Lemma 3.3}
Let $\Bbb{L}$ be a semisimple Lie algebra, H a fixed maximal torus, $\Phi$ the root system of $\Bbb{L}$ corresponding to $H$,
$\Delta=\{\alpha_{1},\alpha_{2},\cdots,\alpha_{l}\}$ a basis of $\Phi$,
and $L$ a minimal Q-graded subalgebra of $\Bbb{L}$ with the following decomposition:
$$ L=H\oplus\bigoplus_{\alpha\in\Psi}L_{\alpha},\,\, \Psi\subseteq\Phi .$$

If $D$ is a derivation of $L$, and it satisfies the condition: $D(H)\subseteq [L,L]$,
then $D(L_{\alpha})\subseteq L_{\alpha}, \forall \alpha\in \Psi$.
\endproclaim

\demo{Proof} Set $I=[L,L]=\oplus_{\alpha\in \Psi}L_{\alpha}$, and $\Psi=\{\beta_1,\cdots,\beta_m\}$, here $m=dim I$.
For any root $\alpha\in \Psi$, and non-zero vector $x_{\alpha}\in L_{\alpha}$, we have
$$ D(x_{\alpha})=y_{\beta_1}+y_{\beta_2}+\cdots +y_{\beta_m}\in I, \,\, y_{\beta_i}\in L_{\beta_i}, \,\,1\leq i\leq m.$$

Since $I$ is an abelian subalgebra by Theorem 2.4, from the assumptions of the lemma, the following identities hold:
$$ \alpha(h)D(x_{\alpha})=D[h,x_{\alpha}]=[D(h),x_{\alpha}]+[h,D(x_{\alpha})]=[h, D(x_{\alpha})],$$
$$ \alpha(h)(y_{\beta_1}+y_{\beta_2}+\cdots +y_{\beta_m})=\beta_1(h)y_{\beta_1}+\beta_2(h)y_{\beta_2}+\cdots +\beta_m(h)y_{\beta_m},$$
$$ \alpha(h)y_{\beta_i}=\beta_i(h)y_{\beta_i}, \,\, 1\leq i\leq m, \,\, \forall h\in H.$$
Hence we have the equality $\beta_i=\alpha$ whenever $y_{\beta_i}\neq 0$, from which the result of the lemma follows immediately.
\hfill $ $ \qed
\enddemo
\proclaim{Lemma 3.4} In the notations of Lemma 3.3 and its proof, if the derivation $D$ of $L$ satisfies the stronger
condition: $D(h)\in [h, L], \forall h\in H$, then there exists some element $z\in L$ such that the following relations hold
$$ (D+ad(z))(H)=0, \,\, (D+ad(z))(L_{\alpha})\subseteq L_{\alpha}, \,\, \forall \alpha\in \Psi .$$
\endproclaim
\demo{Proof} According to the result of Lemma 3.3, the derivation $D$ acts on the abelian subalgebra $I=[L,L]$ as follows:
$$ D(x_{\alpha})=a_{\alpha}x_{\alpha}, \,\, \forall x_{\alpha}\in L_{\alpha}, \forall \alpha\in \Psi,$$
here $a_{\alpha}\in \Bbb{F}, \forall \alpha\in \Psi$. Furthermore, without loss of generality, we may assume that the derivation $D$
acts on $H$ like this: for $h\in H$, we have
$$ D(h)=y_{\beta_1}+y_{\beta_2}+\cdots +y_{\beta_m}, \,\, y_{\beta_i}\in L_{\beta_i}, \,\, \beta_i\in \Psi,$$
here $\beta_{i}(h)\neq 0$ whenever $y_{\beta_i}\neq 0$ for any $i\in \{1,2,\cdots, m\}$.

Now choose a basis $\{h_1,h_2,\cdots ,h_l\}$ of $H$, we will adjust the actions of $D$ on these basis elements,
the process is accomplished in several steps below.

Step 1: Suppose that $D(h_1)\neq 0$ (otherwise go to the next step), and has the same form as the one in the previous paragraph:
$$ D(h_1)=y_{\beta_1}+y_{\beta_2}+\cdots +y_{\beta_m}\in I.$$
Set $z_1=\varepsilon_1\beta_1(h_1)^{-1}y_{\beta_1}+\cdots +\varepsilon_m\beta_{m}(h_1)^{-1}y_{\beta_m}\in I$,
where $\varepsilon_i=1$ if $\beta_i(h_1)\neq 0$ and $\varepsilon_i=0$ otherwise.
Define $D_1=D+ad(z_1)$, it is easy to see that
$$ D_1(x_\alpha)=D(x_{\alpha})=a_{\alpha}x_{\alpha}, \,\, \forall x_{\alpha}\in L_{\alpha}, \forall \alpha\in \Psi,$$
$$ D_1(h_1)=D(h_1)+[z_1, h_1]=0.$$

Step 2: If $D_1(h_2)\neq 0$, similar to step 1, we find some element $z_2\in I$ and define the derivation $D_2=D_1+ad(z_2)$
which obeys the rules:
$$ D_2(x_{\alpha})=D_1(x_{\alpha})=D(x_{\alpha}),\forall x_{\alpha}\in L_{\alpha}, \forall \alpha\in \Psi,$$
$$ D_2(h_2)=D_1(h_2)+[z_2, h_2]=0.$$
Before going to the next step, we must verify that $D_2(h_1)=0$, and this can be seen from the following identities
$$ \eqalign{ D_1(h_2) & =w_{\beta_1}+\cdots +w_{\beta_m}, w_{\beta_i}\in L_{\beta_i}, \beta_i\in \Psi, \cr
                 z_2 & =\varepsilon_1\beta_1(h_2)^{-1}w_{\beta_1}+\cdots +\varepsilon_m\beta_m(h_2)^{-1}w_{\beta_m},\cr
                   0 & =D_1[h_1,h_2]=[D_1(h_1),h_2]+[h_1,D_1(h_2)]=[h_1, D_1(h_2)]\cr
                     & =\beta_1(h_1)w_{\beta_1}+\cdots +\beta_m(h_1)w_{\beta_m}. \cr
               D_2(h_1)& =D_1(h_1)+[z_2,h_1]=[z_2,h_1]\cr
                     & =-\beta_1(h_1)\varepsilon_1\beta_1(h_2)^{-1}w_{\beta_1}-\cdots -\beta_m(h_1)\varepsilon_m\beta_m(h_2)^{-1}w_{\beta_m}\cr
                     & =0,\cr}$$
here $\varepsilon_i=1$ or $0$, $1\leq i\leq m$, and they are similarly defined as in step 1.

Step 3: Repeating the processes in step 1 or 2 above, we will find elements $z_i\in L$ and derivations $D_i\in Der(L)$ satisfying
$$ D_i(x_{\alpha})=D(x_{\alpha})=a_{\alpha}x_{\alpha},\forall x_{\alpha}\in L_{\alpha}, \forall \alpha\in \Psi, $$
$$ D_i(h_1)=D_i(h_2)=\cdots =D_i(h_i)=0.$$
For $i=l$, we define the element $z$ to be the sum of all $z_i$ that appeared in all these steps,
it is obvious that the element $z$ is the required one.
\hfill $ $ \qed
\enddemo

By this time, we have well prepared to prove the main result of this section.
\proclaim{Theorem 3.5} Let $L$ be a minimal Q-graded subalgebra of a semisimple Lie algebra,
then all almost inner derivations of $L$ are inner: $AID(L)=Inn(L)$.
\endproclaim
\demo{Proof} Given an almost inner derivation $D\in AID(L)$, it clearly satisfies the conditions in Lemma 3.3 and 3.4,
so there exists an element $z\in L$ such that the derivation $\tilde{D}=D+ad(z)$ acts on $L$ as follows
$$ (D+ad(z))(H)=0, \,\, (D+ad(z))(L_{\alpha})\subseteq L_{\alpha}, \,\, \forall \alpha\in \Psi .$$

It suffices to show that the derivation $\tilde{D}$ is inner. Since $\tilde{D}$ satisfies the relations above,
we may choose scalars $a_{\alpha}\in \Bbb{F}$, non-zero elements $x_{\alpha}\in L_{\alpha},  \alpha\in \Psi$,
such that $\tilde{D}$ is determined as follows
$$ \tilde{D}(H)=0, \,\, \tilde{D}(x_{\alpha})=a_{\alpha}x_{\alpha}, \forall \alpha\in \Psi.$$

As $\tilde{D}$ is almost inner, according to the definition of almost inner derivations, there exists an element $h\in H$ satisfying
the following identity:
$$ \tilde{D}(\sum_{\alpha\in \Psi}x_{\alpha})=[h, \sum_{\alpha\in \Psi}x_{\alpha}]=\sum_{\alpha\in \Psi}\alpha(h)x_{\alpha} $$
from which we deduce directly that $a_{\alpha}=\alpha(h), \forall \alpha\in \Psi$, and $\tilde{D}$ is the inner derivation $ad(h)$,
this implies that $D$ is an inner derivation.
\hfill $ $ \qed
\enddemo

Remark: Under the conditions of Lemma 3.3 and 3.4 (weaker than that of almost inner derivations),
the modified derivation $\tilde{D}$ appeared in the proof of Theorem 3.5 annihilates the maximal torus $H$ and acts as scalars on root spaces.
Conversely, any endomorphism $\delta$ of $L$ thus defined must be a derivation of $L$.
To make further decisions like whether or not this derivation $\delta$ is inner, we obtain two different results in the two cases below.

a) $\dim [L,L]=m=l=\dim H$, all such $\delta$ are inner derivations;

b) $\dim [L,L]=m>l=\dim H$, there exists non-inner $\delta$ as described above.

To understand better these two results, we suppose that the derivation $\delta$ is given precisely as follows
$$ \delta(H)=0, \,\, \delta(x_{\alpha})=a_{\alpha}x_{\alpha}, \,\,x_{\alpha}\in L_{\alpha},\,\, \alpha\in \Psi,$$
here the scalars $a_{\alpha}$ are only depended on $\alpha\in \Psi$. If $\delta$ is inner, then there exists some element $z\in L$
such that $0=\delta(h)=[z,h] \,(\forall h\in H)$ from which we deduce that $z\in H$.
Hence, we have $\delta(x_{\alpha})=[z, x_{\alpha}]=\alpha(z)x_{\alpha}, \forall \alpha\in \Psi$.
Therefore, the derivation $\delta$ is inner if and only if the following statement holds:
$$ (\beta_1(z), \beta_2(z), \cdots, \beta_m(z))=(a_{\beta_1}, a_{\beta_2},\cdots, a_{\beta_m}), \,\, \exists z\in H.$$

Since all the vectors in the form of left side of the identity constitute a vector space with dimension $l$,
but all the possible vectors in the right side form a vector space with dimension $m$,
we get the two results a) and b) immediately.

 \vskip 5pt

\smallskip\bigskip
\subhead 4. \ \  Derivations of the loop algebra of a minimal Q-graded subalgebra
\endsubhead
\medskip

Recall that any Lie algebra $L$ over the field $\Bbb{F}$ with a fixed invariant symmetric bilinear form $(,)$
has an affinization $\tilde{L}$ defined by the following
$$ \tilde{L}=L\otimes \Bbb{F}[t,t^{-1}]\oplus \Bbb{F}K,$$
where $\Bbb{F}[t,t^{-1}]$ is the commutative associative algebra consisting of Laurent polynomials over $\Bbb{F}$,
and $K$ is a central element of $\tilde{L}$. The bracket of $\tilde{L}$ is determined by the formula
$$ [x\otimes t^m, y\otimes t^n]=[x,y]\otimes t^{m+n}+m(x,y)\delta_{m+n,0}K $$
for any elements $x,y\in L$ and integers $m,n\in \Bbb{Z}$.

We will investigate and determine almost inner derivations of $\tilde{L}$ when $L$ is
a minimal Q-graded subalgebra of a semisimple Lie algebra.
We first consider the situation in loop algebra which is the main part of the
affinization, neglecting the central element $K$, and need some notations and basic results from [2,4].

\proclaim{Lemma 4.1} Let $L$ be a minimal Q-graded subalgebra of a semisimple Lie algebra, and set $S=\Bbb{F}[t,t^{-1}]$,
then the following identity holds:
$$ Der(L\otimes S)=D_S(L\otimes S)\oplus D_{L\otimes 1}(L\otimes S),$$
where the subspace $D_S(L\otimes S)$ is defined as
$$ \{d\in Der(L\otimes S); d(x\otimes s_1s_2)=s_1d(x\otimes s_2)=d(x\otimes s_1)s_2, \forall x\in L,\forall s_1,s_2\in S\},$$
and $D_{L\otimes 1}(L\otimes S)$ denote the subspace $ \{\delta\in Der(L\otimes S); \delta(L\otimes 1)=0\}$.
\endproclaim
\demo{Proof} This is a special case of Lemma 2.3 in [2] (Note that the loop algebra $L\otimes S$ can be regarded naturally
as a $S$-$S$-bimodule).
\hfill $ $ \qed
\enddemo

\proclaim{Lemma 4.2} With notions as above, the subspace $D_S(L\otimes S)$ can be identified with the tensor product $Der(L)\otimes S$.
Moreover, any almost inner derivation in the subspace $D_S(L\otimes S)$ is an inner derivation.

\endproclaim
\demo{Proof} It is obvious that a derivation $d\in D_S(L\otimes S)$ is determined by its action on the subalgebra $L\otimes 1$.
Since $\dim(L)<\infty $, for any element $x\in L$, we can write
$$ d(x\otimes 1)=\sum_{i=r_1}^{r_2} D_i(x)\otimes t^i, $$
where $D_i: L\rightarrow L$ is a derivation of $L$, $r_1\leq i\leq r_2$,  and $r_1,r_2$ are integers
which are depended only on the dimension of $L$. Therefore, we have
$$ d(x\otimes s)=d(x\otimes 1)s=\sum_{i=r_1}^{r_2} D_i(x)\otimes t^i s=(\sum_{i=r_1}^{r_2} D_i\otimes t^i)(x\otimes s) $$
for any $x\in L, s\in S$, here the symbol $t^i$ is also used to denote the related left multiplication of $S$.
This proves the first statement of the lemma.

Now let $d$ be an almost inner derivation in the subspace $D_S(L\otimes S)$, there must exists some element $\sum_i y_i\otimes t^i\in L\otimes S$ satisfying
$$ d(x\otimes 1)=[x\otimes 1,\sum_i y_i\otimes t^i]=\sum_i [x,y_i]\otimes t^i $$
from which it follows that $D_i(x)=[x,y_i]$ and $D_i$ is an almost inner derivation of the Lie algebra $L$,
hence $D_i$ is an inner derivation of $L$ by Theorem 3.5, $r_1\leq i\leq r_2$. Suppose that $D_i=ad(-y_i), r_1\leq i\leq r_2$,
we deduce that
$$ d(x\otimes 1)=ad(-\sum_{i=r_1}^{r_2} y_i\otimes t^i)(x\otimes 1), \,\, \forall x\in L.$$

Set $w=-\sum_{i=r_1}^{r_2} y_i\otimes t^i$ and take a generic element $\sum_k x_k\otimes t^k\in L\otimes S$, we have
$$ d(\sum_k x_k\otimes t^k)=\sum_k d(x_k\otimes 1)t^k $$
$$ =\sum_k [w,x_k\otimes 1]t^k=ad(w)(\sum_k x_k\otimes t^k).$$
Therefore the derivation $d$ is an inner derivation as required.
\hfill $ $ \qed
\enddemo

\proclaim{Definition 4.3} For any non-associative algebra $A$ over the field $\Bbb{F}$, define the centroid $Cent(A)$ of $A$ as follows
$$ Cent(A)=\{\varphi\in End A; \varphi(xy)=\varphi(x)y=x\varphi(y), x,y\in A\}.$$
Obviously, $Cent(A)$ is an associative subalgebra with unit of $End A$. when $A=L$ is a Lie algebra over $\Bbb{F}$,
the condition for $\varphi$ is written as usual
$$\varphi[x,y]=[\varphi(x),y]=[x,\varphi(y)], x,y\in L. $$
\endproclaim

The following lemma presents a similar result to Lemma 1 (Ch.X, \S 1) in [11] which is not applicable here,
since the algebra $L$ in our case is not perfect, so we will give a direct proof of it.

\proclaim{Lemma 4.4} Suppose that $L$ is a minimal Q-graded subalgebra of a semisimple Lie algebra
with the additional assumption that $\dim H=\dim [L,L]$, then $Cent(L)$ is a commutative algebra.
\endproclaim
\demo{Proof} Decompose $L$ as a direct sum of root subspaces as before
$$ L=H\oplus\bigoplus_{\alpha\in \Psi}L_{\alpha},$$
choose a basis $\{h_1, \cdots, h_l\}$ of $H$ which is dual to the basis $\Psi=\{\beta_1,\cdots , \beta_l\}$ of
the dual space $H^*$, and fix nonzero vectors $x_i\in L_{\beta_i}, 1\leq i\leq l$.

Let $\varphi\in Cent(L)$, then $\varphi(x_j)\in [L,L]$ since $x_j\in [L, L]$,
hence we may suppose that $\varphi(x_j)=\sum_k a_{kj}x_k, a_{kj}\in \Bbb{F}$, and we have
$$ \varphi[h_i, x_j]=[h_i, \varphi(x_j)]=\sum_k a_{kj}\beta_k(h_i)x_k=a_{ij}x_i$$
and
$$ \varphi[h_i, x_j]=\varphi(\beta_j(h_i)x_j)=\delta_{j,i}\varphi(x_j)$$
from which we deduce that $a_{ij}=0$ whenever $i\neq j$.

Since $0=\varphi[h,h']=[\varphi(h),h']$ for any $h,h'\in H$, we see that $\varphi(h)\in H, \forall h\in H$.
Set $\varphi(h_i)=\sum_k b_{ki}h_k, b_{ki}\in \Bbb{F}$, then similarly we have the following
$$ \varphi[h_i, x_j]=[\varphi(h_i), x_j]=[\sum_k b_{ki}h_k, x_j]=b_{ji}x_j$$
from which we deduce that $b_{ij}=0$ whenever $i\neq j$, and $a_{ii}=b_{ii}, \forall i$.
So the linear map $\varphi$ corresponds to a diagonal matrix with respect to the prescribed basis of $L$,
and the associative algebra $Cent(L)$ is commutative.
\hfill $ $ \qed
\enddemo

The following result concerning centroids of algebras is fundamental to handle derivations of loop algebras,
there are many similar results in literature. Lemma 1.2 in [2] got the same result with the first tensor factor
being any perfect or unital algebra, etc.  Remark 2.22 in [4] mentioned a general related result, and detailed discussions
appeared in [1].  Here, we deal with only loop algebras of minimal Q-graded subalgebras (not perfect) of semisimple Lie algebras,
a direct calculation can be given to get the result.

\proclaim{Lemma 4.5} Let $L$ be given as in Lemma 4.4, then the centroid of the loop algebra $L\otimes S$ is
an associative commutative algebra, and it can be identified with the tensor product algebra $Cent(L)\otimes S$
$$ Cent(L\otimes S)\approx Cent(L)\otimes S. $$
\endproclaim
\demo{Proof} Choose a basis $\{h_i\otimes t^j, x_i\otimes t^j; 1\leq i\leq l, j\in \Bbb{Z}\}$ of $L\otimes S$,
here $\{h_1,\cdots , h_l\}$ is a basis of $H$, $\Psi=\{\beta_1,\cdots, \beta_l\}$ is the dual basis of $H^*$,
and $\{x_1,\cdots, x_l\}$ is a basis of $[L,L]$ corresponding to $\beta_1,\cdots, \beta_l $, so $L$ has the following decomposition
$$ L=H\oplus L_{\beta_1}\oplus \cdots \oplus L_{\beta_l}.$$

Recall that $S=\Bbb{F}[t,t^{-1}]$ as in Lemma 4.1, and if $\varphi\in Cent(L\otimes S)$, then we have
$$ [\varphi(h_i\otimes t^j), h'\otimes 1]=\varphi[h_i\otimes t^j, h'\otimes 1]=0, \forall h'\in H,$$
hence $\varphi(h_i\otimes t^j)\in H\otimes S, \forall i,j$.
Also, $x_i\in [L,L]$ implies that $\varphi(x_i\otimes t^j)\in [L,L]\otimes S$ and so we can write
$$ \varphi(h_i\otimes t^j)=\sum_k h_k\otimes f_{ij}^k(t), \,\,\varphi(x_i\otimes t^j)=\sum_k x_k\otimes g_{ij}^k(t),\,\,\forall i,j,$$
where $f_{ij}^k(t), g_{ij}^k(t)\in S$ are Laurent polynomials, \,\, $\forall i,j,k$.

According to the definition of $\varphi$, the following identities hold
(note that the bases $\{h_1,\cdots,h_l\}$ and $\{\beta_1,\cdots,\beta_l\}$ are dual to each other):
$$ \varphi[h_i\otimes t^j, x_m\otimes t^n]=[\varphi(h_i\otimes t^j), x_m\otimes t^n]=[h_i\otimes t^j, \varphi(x_m\otimes t^n)],$$
$$ [\varphi(h_i\otimes t^j), x_m\otimes t^n]=[\sum_k h_k\otimes f_{ij}^k(t), x_m\otimes t^n]=x_m\otimes f_{ij}^m(t) t^n, $$
$$ [h_i\otimes t^j, \varphi(x_m\otimes t^n)]=[h_i\otimes t^j, \sum_k x_k\otimes g_{mn}^k(t)]=x_i\otimes g_{mn}^i(t) t^j $$
from which we deduce the following relations
$$ f_{ij}^m(t)=g_{ij}^m(t)=0, \,\,i\neq m; \,\,g_{m,j+n}^m(t)=f_{mj}^m(t)t^n, \,\, f_{ij}^i(t)=g_{ij}^i(t),$$
and thus we can adjust the action of the map $\varphi\in Cent(L\otimes S)$ simply as follows
$$ \varphi(h_i\otimes t^j)=h_i\otimes f_{ij}(t), \,\,\varphi(x_i\otimes t^j)=x_i\otimes f_{ij}(t),$$
$$ f_{i,j+k}(t)=f_{ij}(t) t^k=f_{i0}(t) t^{k+j}, \,\,  1\leq i\leq l, \,\,j,k\in \Bbb{Z}.$$
From these formulas we conclude that the algebra $Cent(L\otimes S)$ is commutative.

Conversely, given arbitrary Laurent polynomials $f_{i0}\in \Bbb{F}[t,t^{-1}] \,(1\leq i\leq l)$ and set $f_{ij}(t)=f_{i0}(t)t^j$,
then we define a linear map $\varphi: L\otimes S\rightarrow L\otimes S$ as above,
it can be checked that $\varphi\in Cent(L\otimes S)$.

Finally, we show that $Cent(L\otimes S)$ is isomorphic to $Cent(L)\otimes S$,
and an isomorphism between the two algebras is given below:
$$\sigma: Cent(L)\otimes S\rightarrow Cent(L\otimes S), \lambda\otimes f(t)\mapsto \lambda\otimes L_{f(t)},$$
here the map $L_{f(t)}: S\rightarrow S$ is defined by left multiplication by $f(t)$, $\lambda\in Cent(L)$
is determined by the following \,(Lemma 4.4)
$$ \lambda(h_i)=a_ih_i, \lambda(x_i)=a_ix_i, \,\,a_i\in \Bbb{F}, \, 1\leq i\leq l.$$

To show that the map $\sigma$ is surjective, we define $\lambda_i\in Cent(L)$ by the formula:
$$ \lambda_i(h_j)=\delta_{ij}h_j, \lambda_i(x_j)=\delta_{ij}x_j, \forall i,j,$$
then $\{\lambda_1,\cdots,\lambda_l\}$ is a basis of $Cent(L)$. For any given scalars $b_{ij}\in \Bbb{F}$, define a linear map $\varphi$ by the rule below:
$$ \varphi(h_m\otimes t^n)=\sigma(\sum_{i,j}b_{ij}\lambda_i\otimes t^j)(h_m\otimes t^n)=h_m\otimes \sum_jb_{mj}t^{j+n},$$
$$ \varphi(x_m\otimes t^n)=\sigma(\sum_{i,j}b_{ij}\lambda_i\otimes t^j)(x_m\otimes t^n)=x_m\otimes \sum_jb_{mj}t^{j+n} $$
from which we deduce that the homomorphism $\sigma$ is surjective. Also, it is easy to show that the map $\sigma$ is injective
and the proof of the lemma is finished.
\hfill $ $ \qed
\enddemo

Recall that a derivation from an associative or Lie algebra $A$ to its module $M$ is a linear map
$\delta: A\rightarrow M$ satisfying the following conditions:
$$ \delta(ab)=a\delta(b)-b\delta(a), \,\, \forall a,b\in A.$$

Particularly, derivations from a subalgebra $B$ of the given algebra $A$ to $A$ itself
is defined as derivations from algebras to their modules. The space of all the derivations from $B$ to $A$ is denoted by $D(B,A)$.

\proclaim{Lemma 4.6} Let $L$ be given as in Lemma 4.4 and $S=\Bbb{F}[t,t^{-1}]$, then there is an isomorphism of vector spaces
$$ \tau: D(1\otimes S, Cent(L)\otimes S)\rightarrow D_{L\otimes 1}(L\otimes S), \,\, d\mapsto \tau(d),$$
here $\tau(d)(x\otimes f)=\sigma(d(1\otimes f))(x\otimes 1), \forall x\in L, \forall f\in S$,
and $\sigma $ is the isomorphism given in the proof of Lemma 4.5.
\endproclaim
\demo{Proof} The proof of this lemma can be copied almost word for word from the proof of Lemma 2.4 in [2] using Lemma 4.5,
though the Lie algebra $L$ is not perfect here in our case.
\hfill $ $ \qed
\enddemo

\proclaim{Lemma 4.7} Notations as in Lemma 4.5, then any derivation $D\in D_{L\otimes 1}(L\otimes S)$ can be determined as follows:
$$ D(h_i\otimes t^j)=h_i\otimes jt^{j-1}f_i(t), 1\leq i\leq l, \,\, j\in \Bbb{Z},$$
$$ D(x_i\otimes t^j)=x_i\otimes jt^{j-1}f_i(t), 1\leq i\leq l, \,\, j\in \Bbb{Z},$$
where $f_1(t),\cdots, f_l(t)\in S$ are Laurent polynomials.

In particular, if $D\in D_{L\otimes 1}(L\otimes S)$ is an almost inner derivation of the loop algebra $L\otimes S$,
then it is zero.
\endproclaim
\demo{Proof} By Lemma 4.6, there exists some derivation $d\in D(1\otimes S, Cent(L)\otimes S)$ such that $D=\tau(d)$.
The derivation $d$ is determined by its action on the generator $1\otimes t$ which can be written as
$$ d(1\otimes t)=\sum_{i=1}^l\lambda_i\otimes f_i(t),$$
where $\lambda_1,\cdots, \lambda_l$ is a basis of $Cent(L)$ as defined in Lemma 4.5, $f_1(t),\cdots, f_l(t)$ are Laurent polynomials in $S$,
and we obtain the values of the derivation $D$ in basis elements as required:
$$ \eqalign{ & D(h_i\otimes t^j)=\sigma(d(1\otimes t^j))(h_i\otimes 1) \cr
             & =j(1\otimes t^{j-1})(\sum_{k=1}^l\lambda_k\otimes f_k(t))(h_i\otimes 1)=h_i\otimes jt^{j-1}f_i(t),\cr
             & D(x_i\otimes t^j)=\sigma(d(1\otimes t^j))(x_i\otimes 1) \cr
             & =j(1\otimes t^{j-1})(\sum_{k=1}^l\lambda_k\otimes f_k(t))(x_i\otimes 1)=x_i\otimes jt^{j-1}f_i(t).\cr}$$

If $D\in D_{L\otimes 1}(L\otimes S)$ is an almost inner derivation of $L\otimes S$,
then by the formula above and $D(h_i\otimes t^j)\in [L,L]\otimes S$, we have $f_i(t)=0$, $\forall i$, hence $D=0$.
\hfill $ $ \qed
\enddemo

\proclaim{Theorem 4.8} Suppose that $L$ is a minimal Q-graded subalgebra of a semisimple Lie algebra
with the additional assumption that $\dim H=\dim [L,L]$, then the derivation algebra $Der(L\otimes S)$ has a decomposition:
$$ Der(L\otimes S)=Der(L)\otimes S\oplus Cent(L)\otimes Der(S).$$
Moreover, any almost inner derivation of the loop algebra $L\otimes S$ is inner.
\endproclaim
\demo{Proof} According to Lemma 4.1 and 4.2,
to obtain the decomposition formula in the theorem, it suffices to verify that $D_{L\otimes 1}(L\otimes S)=Cent(L)\otimes Der(S)$.
Using notations in Lemma 4.7, an element $D\in D_{L\otimes 1}(L\otimes S)$ is determined by Laurent polynomials
$f_1(t),\cdots, f_l(t)\in S$ such that
$$ D(h_i\otimes t^j)=h_i\otimes jt^{j-1}f_i(t)$$
$$=\sum_k \lambda_k(h_i)\otimes d_k(t^j)=(\sum_k \lambda_k\otimes d_k)(h_i\otimes t^j), $$
$$ D(x_i\otimes t^j)=x_i\otimes jt^{j-1}f_i(t)$$
$$=\sum_k \lambda_k(x_i)\otimes d_k(t^j)=(\sum_k \lambda_k\otimes d_k)(x_i\otimes t^j), $$
here $d_k\in Der(S)$ is determined by $d_k(t)=f_k(t)$, $1\leq k\leq l$ and $1\leq i\leq l,\, j\in \Bbb{Z}$.
Hence we see that $D\in Cent(L)\otimes Der(S)$ and $D_{L\otimes 1}(L\otimes S)\subset Cent(L)\otimes Der(S)$.
The opposite direction of this relation can be verified easily.

By Lemma 4.1, any derivation $D$ of the loop algebra $L\otimes S$ can be written as a sum,
$$ D=d+(D-d),$$
where $d\in D_S(L\otimes S)$ is defined by the formula: $d(x\otimes f)=D(x\otimes 1)f, x\in L, f\in S$,
 and $D-d\in D_{L\otimes 1}(L\otimes S)$.

If $D$ is an almost inner derivation, then both the derivations $d$ and $D-d$ are almost inner derivations.
It suffices to show that $d$ is almost inner, this can be seen from the following identities:
$$ \eqalign{ & d(\sum_i h_i\otimes f_i)=\sum_i D(h_i\otimes 1)f_i =\sum_i [h_i\otimes 1, \sum_k x_k\otimes g_{ik}]f_i \cr
             & =\sum_i x_i\otimes g_{ii}f_i =[\sum_i h_i\otimes f_i, \sum_k x_k\otimes g_{kk}],\cr
             & d(\sum_i x_i\otimes \varphi_i)=\sum_i D(x_i\otimes 1)\varphi_i =\sum_i [x_i\otimes 1, \sum_k h_k\otimes \psi_{ik}]\varphi_i \cr
             & =-\sum_i x_i\otimes \psi_{ii}\varphi_i =[\sum_i x_i\otimes \varphi_i, \sum_k h_k\otimes \psi_{kk}], \cr
             & d(\sum_i h_i\otimes f_i+\sum_i x_i\otimes \varphi_i) \cr
             & =[\sum_i h_i\otimes f_i, \sum_k x_k\otimes g_{kk}]+[\sum_i x_i\otimes \varphi_i, \sum_k h_k\otimes \psi_{kk}]\cr
             & =[\sum_i h_i\otimes f_i+\sum_i x_i\otimes \varphi_i, \sum_k x_k\otimes g_{kk}+\sum_k h_k\otimes \psi_{kk}],\cr} $$
where $\{h_1,\cdots, h_l, x_1,\cdots, x_l\}$ is the basis of $L$ as defined in Lemma 4.5 and 4.7, and $f_i, g_{ik}, \varphi_i$, $\psi_{ik}$ are
Laurent polynomials in $S$ \,$(1\leq i,k\leq l)$.

We conclude that $D$ is inner from Lemma 4.2 and Lemma 4.7.
\hfill $ $ \qed
\enddemo

\vskip 5pt

\smallskip\bigskip
\subhead 5. \ \  Derivations of the affinization of a minimal Q-graded subalgebra
\endsubhead
\medskip

We are ready to investigate derivations of the affinization of a minimal Q-graded subalgebra of a semisimple Lie algebra $\Bbb{L}$.
So let $L$ be a minimal Q-graded subalgebra with $\dim L=2\dim [L,L]$ as in Lemma 4.4,
the affinization of $L$ is the Lie algebra $\tilde{L}$ which is defined as follows:
$$ \tilde{L}=L\otimes S \oplus \Bbb{F}K, $$
where $S=\Bbb{F}[t,t^{-1}]$ is the algebra of Laurent polynomials over $\Bbb{F}$,
and $K$ is a central element of $\tilde L$. The bracket of $\tilde L$ is given below:
$$ [x\otimes t^m, y\otimes t^n]=[x,y]\otimes t^{m+n}+m(x,y)\delta_{m+n,0}K $$
for any elements $x,y\in L$ and integers $m,n\in \Bbb{Z}$.

Note that the bilinear form $(,)$ is the restriction of the Killing form of $\Bbb{L}$ to its subalgebra $L$,
and $\Bbb{F}K$ is actually the center of the Lie algebra $\tilde{L}$. Note also that any derivation of $\tilde{L}$
maps the center $\Bbb{F}K$ into itself.

\proclaim{Lemma 5.1} Let $D$ be an almost inner derivation of the Lie algebra $\tilde{L}$, then there exists an element $y\in \tilde{L}$,
such that the following relations hold:
$$ (D-ad(y))(I\otimes S)=0, \,\, (D-ad(y))(K)=0, $$
$$ (D-ad(y))(H\otimes S)\subset \Bbb{F}K. $$
\endproclaim
\demo{Proof} Since $D$ is almost inner, we have $D(K)\in [K, \tilde{L}]$ and $D(K)=0$.
Thus the derivation $D$ induces a map $\bar{D}$ as follows
$$ \bar{D}: \tilde{L}/\Bbb{F}K\rightarrow \tilde{L}/\Bbb{F}K, \,\, \bar{w}\mapsto \overline{D(w)}.$$
It can be checked that the map $\bar{D}$ is an almost inner derivation of the quotient algebra $\tilde{L}/\Bbb{F}K$,
which is isomorphic to the loop algebra $L\otimes S$.

From Theorem 4.8, we deduce that there exists some element $y\in \tilde{L}$ such that $\bar{D}(\bar{w})=ad(\bar{y})(\bar{w}), \forall w\in \tilde{L}$,
and therefore we have $(D-ad(y))(\tilde{L})\subset \Bbb{F}K$.
In particular, we obtain that $(D-ad(y))[\tilde{L}, \tilde{L}]=0 $ and all the relations of the lemma hold for this element $y$.
\hfill $ $ \qed
\enddemo

\proclaim{Theorem 5.2} Let $\tilde{L}$ be the affinization of the minimal Q-graded subalgebra $L$ of a semisimple Lie algebra
with $\dim{L}=2\dim [L,L]$, then $D\in End \tilde{L}$ is an almost inner derivation if and only if
there exists an element $y\in \tilde{L}$ such that the corresponding relations in Lemma 5.1 hold.
\demo{Proof} By Lemma 5.1, it suffices to show that the linear map $D$ is an almost inner derivation if it satisfies the following conditions:
$$ D(I\otimes S)=0, \,\, D(K)=0, \,\, D(H\otimes S)\subset \Bbb{F}K. $$

It is clear that $D$ is a derivation if these conditions are satisfied. To show that $D$ is almost inner, we do
some concrete calculations below.
Fix a basis $\{h_1,\cdots,h_l, x_1,\cdots, x_l\}$ of the Lie algebra $L$ as in Lemma 4.7,
and choose a dual basis $\{h_1',\cdots ,h_l'\}$ of $\{h_1,\cdots, h_l\}$ in $H$ with respect to the given bilinear form $(,)$,
so we have $(h_i, h_j')=\delta_{ij}, 1\leq i,j\leq l$.

Note that $h_i'=t_{\beta_i}, \beta_i(h_j)=(t_{\beta_i}, h_j)=\delta_{ij}, \forall i,j$, and $\Psi=\{\beta_1,\cdots, \beta_l\}$ as in Lemma 4.5.
Moreover, we have (see [10] for some notations used here):
$$ \beta_i(h_j')=(t_{\beta_i}, t_{\beta_j})=(\beta_i, \beta_j), \,\, 1\leq i,j\leq l.$$

Set $D_{ij}(h_m\otimes t^n)=\delta_{im}\delta_{jn}K \,(1\leq m\leq l, n\in \Bbb{Z})$, $D_{ij}(I\otimes S)=D_{ij}(K)=0$,
it extends to a unique element $D_{ij}\in End(\tilde{L})$ that satisfies the same relations as $D$ does, and we have
$$ D=\sum_{m,n}a_{mn}D_{mn}, \,\, D(h_m\otimes t^n)=a_{mn}K, \,\, \forall m,n. $$
Though the sum is infinite, it is well defined.
We will show that each $D_{ij}$ is almost inner, this suffices for the given derivation $D$ to be almost inner.

Claim: $D_{ij}$ is an almost inner derivation.

Without loss of generality, we may choose any element with no
central part: $X=\sum_{m,n}b_{mn}h_m\otimes t^n +\sum_{u,v}c_{uv}x_u\otimes t^v \in \tilde{L}, b_{mn}, c_{uv}\in \Bbb{F}$, and define
$$ A=\{m\in \{1,\cdots, l\}; \exists n, s.t.,  b_{mn}\neq 0\}, \,\, B=\{1,\cdots, l\}\backslash A.$$
Let $Y=j^{-1}h_i'\otimes t^{-j}-\sum_{k\in B}d_kh_k'\otimes t^{-j}+\sum_{p,q}e_{pq}x_p\otimes t^q\in \tilde{L}$,
here all the coefficients $d_k$ and $e_{pq}$ are to be fixed later to satisfy: $D_{ij}(X)=[X,Y]$.  We calculate the brackets as follows
$$ \eqalign{ & b_{ij}K  =D_{ij}(X)=[X,Y]
    = [\sum_{m\in A, n}b_{mn}h_m\otimes t^n +\sum_{u,v}c_{uv}x_u\otimes t^v, \cr
   &   \,\,\,\,\,\,\,\,\,\, j^{-1}h_i'\otimes t^{-j}-
        \sum_{k\in B}d_kh_k'\otimes t^{-j}+\sum_{p,q}e_{pq}x_p\otimes t^q] \cr
   & = b_{ij}K+\sum_{m\in A, n,q}b_{mn}e_{mq}x_m\otimes t^{n+q}
       -\sum_{u,v}\beta_u(h_i')c_{uv}j^{-1}x_u\otimes t^{v-j} \cr
   & \,\,\,\,\, +\sum_{m\in B,u,v}\beta_u(h_m')c_{uv}d_m x_u\otimes t^{v-j} \cr
   & = b_{ij}K+\sum_{k\in A, n,q}b_{k,n-q}e_{kq}x_k\otimes t^n
       -\sum_{k,n}\beta_k(h_i')c_{k,n+j}j^{-1}x_k\otimes t^n \cr
   & \,\,\,\,\, +\sum_{m\in B,k,n}\beta_k(h_m')c_{k,n+j}d_m x_k\otimes t^n \cr
   & = b_{ij}K+\sum_{k\in A,n}\{\sum_q b_{k,n-q}e_{kq}-\beta_k(h_i')c_{k,n+j}j^{-1}
     +\sum_{m\in B}\beta_k(h_m')c_{k,n+j}d_m\}x_k\otimes t^n \cr
   &  \, \,\,\,-\sum_{k\in B, n}\{\beta_k(h_i')c_{k,n+j}j^{-1}- \sum_{m\in B}\beta_k(h_m')c_{k,n+j}d_m\}x_k\otimes t^n.\cr} $$

For $k\in A$ and any $n$, the following single equation has solutions $e_{kq}$ with arbitrary values of $d_m$,
since its coefficients $b_{k, n-q}$ are not all zero
$$ \sum_q b_{k,n-q}e_{kq}-\beta_k(h_i')c_{k,n+j}j^{-1}+\sum_{m\in B}\beta_k(h_m')c_{k,n+j}d_m=0.$$

For $k\in B$ and any $n$, the following system of equations on unknowns $d_m$
$$ \beta_k(h_i')c_{k,n+j}j^{-1}=\sum_{m\in B}\beta_k(h_m')c_{k,n+j}d_m $$
or equivalently
$$ \beta_k(h_i')j^{-1}=\sum_{m\in B}\beta_k(h_m')d_m, \,\,k\in B$$
has a unique solution, since the matrix $(\beta_k(h_m'))=(\beta_k, \beta_m)$ of its coefficients is non-degenerate.
The proof of the claim is finished.
\hfill $ $ \qed
\enddemo

\proclaim{Proposition 5.3} Notations as in Theorem 5.2, let $D$ be a derivation of $\tilde{L}$ that satisfies the relations:
 $ D(I\otimes S)=0, \,\, D(K)=0, \,\, D(H\otimes S)\subset \Bbb{F}K$.
 If $D$ is an inner derivation, then $D$ is zero.
\endproclaim
\demo{Proof}
Suppose that the element
$Y=\sum_{u,v}d_{uv}h_u\otimes t^v+\sum_{p,q}e_{pq}x_p\otimes t^q\in \tilde{L}$ satisfies $D(X)=[X,Y]$ for any $X\in \tilde{L}$.
Let $D(h_i\otimes t^j)=a_{ij}K$ for some $a_{ij}\in \Bbb{F}$, we have the following identities
$$ \eqalign{  a_{ij}K  & = D(h_i\otimes t^j+x_m\otimes t^n) \cr
  & =[h_i\otimes t^j+x_m\otimes t^n, \sum_{u,v}d_{uv}h_u\otimes t^v+\sum_{p,q}e_{pq}x_p\otimes t^q] \cr
  & =\sum_u d_{u,-j}(h_i,h_u)jK+\sum_qe_{iq}x_i\otimes t^{j+q}-\sum_{v}d_{mv}x_m\otimes t^{n+v}. \cr} $$

Set $m\neq i$, the above equation shows that $e_{iq}=0, \forall q$ and $d_{mv}=0, \forall v$.
Varying $i$ and $m$ if necessary, we deduce that all the elements $e_{pq}$ and $d_{uv}$ are zero,
and this implies that the derivation $D$ is zero.
\hfill $ $ \qed
\enddemo

Conclusion: From the results of this section, we see that the affinization $\tilde{L}$ of a minimal Q-graded
subalgebra $L$ \,(with $\dim L=2\dim[L,L]$) of a semisimple Lie algebra
has many almost inner derivations that are not inner, and
the quotient space $AID(\tilde{L})/Inn(\tilde{L})$ consists of elements of the following form
$$ \sum_{1\leq i\leq l, j\in \Bbb{Z}}a_{ij}D_{ij}+Inn(\tilde{L}),\,\, a_{ij}\in \Bbb{F}, $$
where the infinite sum is well defined and the almost inner derivation $D_{ij}$ is determined by the following rules:
$$ D_{ij}(I\otimes S)=0, \,\, D_{ij}(K)=0, $$
$$ D_{ij}(h_m\otimes t^n)=\delta_{im}\delta_{jn}K,\,\, 1\leq m\leq l, \,\, n\in \Bbb{Z}.$$
Moreover, the infinite subset $\{D_{ij}+Inn(\tilde{L}); 1\leq i\leq l, j\in \Bbb{Z}\}$ of the quotient vector space
$ AID(\tilde{L})/Inn(\tilde{L}) $ is linearly independent\,(Proposition 5.3).

\vskip.3cm \Refs\nofrills{\bf REFERENCES}
\bigskip
\parindent=0.45in

\leftitem{[1]} B.Allison, S.Berman, A.Pianzola, Iterated loop algebras,
 Pacific Journal of Mathematics 227:(1)(2006), 1-41.

\leftitem{[2]} S.Azam, Derivations of tensor product of algebras, Comm. Algebra 36:(3)(2008), 905-927.

\leftitem{[3]} G.Benkart, Derivations and invariant forms of Lie algebras graded by finite root systems,
 Canadian Journal of Mathematics 50(1998), 225-241.

\leftitem{[4]} G.Benkart, E.Neher, The centroid of extended affine and root graded Lie algebras,
 J. of Pure and Applied Algebra 205(2006), 117-145.

\leftitem{[5]} G.Benkart, E.Zelmanov, Lie algebras graded by finite root systems and intersection matrix algebras,
 Invent. Math. 126(1996), 1-45.

\leftitem{[6]} D.Burde, K.Dekimpe, B.Verbeke, Almost inner derivations of Lie algebras,
 J. of Algebras and Its Applications 17(2018), no.11, 26 pages.

\leftitem{[7]} D.Burde, K.Dekimpe, B.Verbeke, Almost inner derivations of Lie algebras II,
ArXiv: 1905.08145(2019).

\leftitem{[8]} D.Burde, K.Dekimpe, B.Verbeke, Almost inner derivations of 2-step nilpotent Lie algebras of genus 2,
ArXiv: 2004.10567v1 [math,RA] 22 Apr 2020.

\leftitem{[9]} C.S.Gordon, E.N.Wilson, Isospectral deformations of compact solvmanifolds,
 J.Differen- tial Geom.19(1984), no.1, 214-256.

\leftitem{[10]} J.Humphreys, Introduction to Lie algebras and representation theory, GTM9, Springer, Fifth edition, 1987.

\leftitem{[11]} N.Jacobson, Lie algebras, Dover, New York, 1979.

\leftitem{[12]} I.Kenji, K.Yoshiyuki, Representation Theory of the Virasoro Algebra,
Springer Monographs in Mathematics, Springer, 2011.

\leftitem{[13]} Dong Liu, Naihong Hu, Leibniz algebras graded by finite root systems,
 Algebra Colloquium 17:(3)(2010), 431-446.

\leftitem{[14]} J. Milnor, Eigenvalues of the Laplace operator on certain manifolds,
Proc.Nat.Acad.Sci. U.S.A.51(1964), 542ff.

\endRefs
\vfill
\enddocument
\end